\documentclass[final,twocolumn,twoside,10pt]{IEEEtran}

\IEEEoverridecommandlockouts
\usepackage{float}

\usepackage{amsmath,mathrsfs, eufrak, mathtools,amsthm}
\usepackage{amsmath,amsfonts,amssymb,amsthm,epsfig,epstopdf,array}
\usepackage{amssymb}
\usepackage{pstricks,pst-plot,psfrag}
\usepackage[all]{xy}
\usepackage{graphicx,subfigure,xspace,bm}
\usepackage[lined,linesnumbered,algoruled,noend]{algorithm2e}

\newtheorem{theorem}{Theorem}[section]
\newtheorem{definition}[theorem]{Definition}

\newtheorem{lemma}[theorem]{Lemma}

\newtheorem{example}[theorem]{Example}

\newcommand{\oprocendsymbol}{\hbox{$\bullet$}}
\newcommand{\oprocend}{\relax\ifmmode\else\unskip\hfill\fi\oprocendsymbol}

\usepackage{tikz}
\usetikzlibrary{external}
\usetikzlibrary{decorations.shapes,arrows,calc,decorations.markings,shapes,decorations.pathreplacing,shapes.arrows}
\usetikzlibrary{positioning}
\tikzset{main node/.style={circle,fill=blue!20,draw,inner sep=1pt},}

\definecolor{BBlue}{cmyk}{.98,0.10,0,.25}

\begin{document}
\title{Optimal Policies for the Sequential Stochastic Threshold Assignment Problem}
\author{Aristomenis Tsopelakos, Sheldon H. Jacobson \thanks{Aristomenis Tsopelakos is with the Coordinated Science Laboratory and the Electrical and Computer Engineering department, University of Illinois at Urbana-Champaign,~\texttt{tsopela2@illinois.edu}.}}

\maketitle

\begin{abstract}
The Stochastic Sequential Threshold Assignment Problem (SSTAP) addresses the optimal assignment of arriving tasks (jobs) to available resources (workers) to maximize a reward function which consists of indicator functions that incorporate threshold constraints. We present an optimal assignment policy for SSTAP, independent of the probability distribution of the job values and of the number of arriving jobs. We show through an example that this type of reward function can model aviation security problems. We analyze the performance limitations of systems that use the SSTAP optimal assignment policy. Finally, we study the multiple levels SSTAP and the SSTAP with uncertainties in workers performance rates.
\end{abstract}

\section{Introduction}\label{sec:intro}
The stochastic sequential assignment problem (SSAP) addresses the assignment of entities (jobs) to available resources (workers) under uncertainties in the parameters. The uncertainties are typically modeled as probability distributions that govern the random parameters of the problem \cite{Bill:79}. SSAP appears naturally in the passenger screening process for aviation security and in the Internet for the optimal assignment of online requests to available servers.

SSAP was introduced in \cite{DLR:72} where an optimal assignment policy is proven for independent and identically distributed (i.i.d.) random job values, based on the partition of the domain of jobs into subintervals. Optimal assignment identifies the subinterval for each job. Kennedy \cite{K:86} presents an updated optimal assignment policy for random job values that are not necessarily i.i.d.. 

Sakaguchi \cite{SAK:83} discusses a generalization of the SSAP for unknown total number of jobs and, Nikolaev and Jacobson \cite{NJ:10} for random total number of jobs. Other variations of the problem include the optimal sequential assignment with random arrival times and reward functions with discount factors \cite{ASC:74}, the optimal policy for SSAP with random deadlines \cite{RIG:87} and the SSAP with uncertainty in the job value distribution \cite{LJ:11}. 

Markov chain is a regularly used mathematical structure for the modeling of sequential random processes. Nakai and Toru \cite{NAK:86} discuss SSAP for partially observed Markov chains. Baharian and Jacobson \cite{BAJA:13} provide a Markov-decision-process approach for the assignment of tasks under a threshold criterion, which minimizes the probability of the total reward failing to a target value. Furthermore, Baharian and Jacobson \cite{BAHJAC:14} obtain stationary policies, which achieve the optimal expected reward per task as the number of tasks approaches infinity, with distributions governed by an ergodic Markov chain.

Apart from the uncertainties in the job values, uncertainties in the performance rates of the workers may occur. In the case of human workers, the performance rates are uncertain due to human errors and fatigue; in the case of machines due to disturbances in measurements and equipment aging. This issue led to the study of SSAP under uncertainties in workers performance rates. The study of optimal assignment policies for SSAP with random worker performance rates \cite{KBKJ:13} and with time-dependent performance rates \cite{BKJ:16} has led to the formulation and solution of the doubly stochastic sequential assignment problem \cite{KHJ:16}.

SSAPs have found application in numerous areas. Nikolaev et al. \cite{NJMcL:07} addresses the sequential stochastic security design problem (SSSDP), which models passenger and carry-on baggage-screening operations in an aviation security system,  to maximize the total security of all passenger-screening decisions over a fixed time period, given passenger risk levels and security device parameters.  

McLay et al. \cite{McLSheAle:09} introduces the Sequential Stochastic Passenger Screening Problem (SSPSP), which allows passengers to be optimally assigned (in real-time) to aviation security resources. Lee et al. \cite{LeeMcLShJ:09} study a real-time sequential binary passenger assignment model as a discrete time difference equation, which is manipulated via nonlinear control techniques. Nikolaev et al. \cite{NLJ:12} address the multistage sequential passenger screening problem (MSPSP) that models passenger and carry-on baggage screening operations in an aviation security system with the capability of dynamically updating the perceived risk of passengers. SSAP also applies to financial problems such as optimal stochastic sequential investment \cite{DLR:75} and investment decisions under uncertainty \cite{Pra:83}. Other applications appear in the fields of computer science for reservation systems \cite{HBV:06} and game theory for online mechanism design \cite{PS:04}.

This paper introduces a version of SSAP based on a new type of reward function defined using indicator functions, which capture threshold constraints. The new version called stochastic sequential threshold assignment problem (SSTAP). We us order-preserving functions in the inequality conditions of the indicator functions and we prove an optimal assignment policy based on a Greedy algorithm which assigns the available resource of smallest value that can satisfy the threshold. We provide an example which illustrates the application of SSTAP in aviation security. A suitable order-preserving function that captures the characteristics of the aviation security problem is used.

Given the performance rates of the workers, the threshold and the order-preserving function, we provide a performance analysis of a system that uses SSTAP optimal assignment policy. We research the maximum and minimum job load that a SSTAP system can service while achieving its maximum reward. In SSTAP the optimal sequential assignment algorithm does not depend on the distributions of the job values. Therefore, we can look for a probability distribution $G_{X}$ for the i.i.d. random job values $\{X_{i},i=1,2,\ldots,n\}$ that maximizes the reward function \eqref{thres_cost_function} for the maximum and minimum job load respectively. The optimal mass probability function is provided and the respective density function is an approximation defined as a summation of Gaussians.

In passenger screening for aviation security, workers are organized in more than one level. We analyze the multiple levels SSTAP, which reflects more accurately the aviation security process, where the workers are partitioned into levels according to their performance rates. Human workers deficiencies or machine workers faults have led to the development of the SSTAP with uncertainties in workers performance rates, which is termed the doubly stochastic sequential threshold assignment problem (DSSTAP). We provide an optimal assignment policy for the DSSTAP.

\section{Preliminaries}

We mention the necessary mathematical background from graph theory and probability theory. We provide an overview of the basic results on SSAP and define the reward function using indicator functions that capture threshold constraints for order-preserving functions.

\subsection{Mathematical background}
We start with basic definitions from graph theory. A graph is a pair $G=(V,E)$, where $V$ is a finite set called the vertex set and $E{\subseteq}V{\times}V$ is the edge set which consists of unordered pairs of vertices. A bipartite graph, is the graph whose set of vertices $V$ is decomposed into two disjoint sets $(S,T)$ such that no two graph vertices within the same set are adjacent. A matching $M$ is a subset of edges $E$ such that each node in $V$ appears in at most one edge in $M$. We provide the notion of the maximum weight bipartite matching used in Section \ref{sec_sstap_rpr}, \cite{DW:2000}.

\begin{definition}(Maximum weight bipartite matching).
	Let us consider a bipartite graph $G=(V,E)$ with bipartition $(S,T)$ and weight function $w:E{\to}\mathbb{R}$. The maximum weight bipartite matching is the matching $M$ that maximizes the matching weight:
	\begin{equation}
	w(M)=\sum_{e{\in}M}w(e)
	\end{equation}
\end{definition}

We also require some basic notions from probability theory. A sequence $(Y_{n}:n{\geq}1)$ of random variables can convergence to a random variable $Y$ in four ways; almost surely, in probability, in the mean square sense and in distribution. We provide the definition of convergence in distribution for a sequence of random variables used in Section \ref{sec_perf_an}, \cite{hb:2015}.

\begin{definition}\label{conv_dist}(Convergence in distribution).
	A sequence $(Y_{n}:n{\geq}1)$ of random variables converges in distribution to random variable $Y$ if $\lim_{n{\to}\infty}F_{Y_n}(x)=F_{Y}(x)$ at all continuity points of $F_{Y}$. Convergence in distribution is denoted by $Y_{n}\xrightarrow{\text{d.}}Y$.
\end{definition}

\subsection{Results on SSAP}
Suppose that there are $n$ workers available to perform $n$ sequentially arriving tasks. The value of task $i$ $\{X_{i},i=1,2,\ldots,n\}$ is a random variable with a known cumulative distribution function $G_{X_{i}}$ and domain $\Omega$. If $\{X_{i},i=1,2,\ldots,n\}$ are independent and identically distributed (i.i.d) random variables, then $G_{X_{i}}{\equiv}G_{\tilde{X}}$. Each worker is characterized by a deterministic performance rate $p_{j}{\in}[0,1]$, $j=1,2,\ldots,n$. Once task $i$ arrives, its value $x_{i}$ is revealed. The goal for SSAP is to invent the optimal assignment policy that maximizes the total expected reward,

\begin{equation}\label{rew_inn_pro}
\max\limits_{{\delta}{\in}{\Delta}_{n}}E\big[ \sum_{i=1}^{n} X_{i}p_{\delta(i)} \big]
\end{equation}
where $\Delta_{n}$ is the set of all permutations of the integers $\{1,2,\ldots,n\}$ and $\delta(i)$ refers to the worker assigned to the $i^{th}$ arriving task. A worker can not be reassigned.

Derman et al. \cite{DLR:72} provide an optimal assignment policy for i.i.d and Kennedy \cite{K:86} for general random variables $\{X_{i},i=1,2,\ldots,n\}$. In both cases, the optimal assignment policy is based on Hardy's Lemma \cite{H:34},

\begin{theorem}(Hardy's Lemma \cite{H:34})\\
	If $a_{1}{\leq}a_{2}{\leq}\ldots{\leq}a_{n}$ and $b_{1}{\leq}b_{2}{\leq}\ldots{\leq}b_{n}$ are sequences of non-negative numbers, then
	\begin{equation}
	\max\limits_{{\delta}{\in}{\Delta}_{n}}\sum_{i=1}^{n}a_{\sigma(i)}b_{i}=\sum_{i=1}^{n}a_{i}b_{i}
	\end{equation}	
\end{theorem}
where $\sigma$ is a permutation of $\{1,2,\ldots,n\}$. 

The optimal assignment policies for i.i.d and for general random job values $\{X_{i},i=1,2,\ldots,n\}$ are given in the following two theorems,

\begin{theorem}\cite{DLR:72}
	For each $n{\geq}1$, there exist numbers $-\infty=a_{0,n}{\leq}a_{1,n}{\leq}a_{2,n}{\leq}\ldots{\leq}a_{n,n}=+\infty$, such that whenever there are $n$ i.i.d jobs values and performance rates $p_{1}{\leq}p_{2}{\leq}\ldots{\leq}p_{n}$ then the optimal choice is to use $p_{i}$ if the random variable of the job is contained in the interval $(a_{i-1,n},a_{i,n}]$. The $a_{i,n}$ are independent of the performance rates but depend on $G_{X}$ according to the following recursive relation
	\begin{equation}
	a_{i,n+1}=\int_{a_{i-1,n}}^{a_{i,n}}zdG_{X}(z) + a_{i-1,n}G(a_{i-1,n}) +a_{i,n}[1-G(a_{i,n})]
	\end{equation}
	for $i=1,2,\ldots,n$, where $-\infty{\cdot}0$ and $\infty{\cdot}0$ are defined to be 0.
\end{theorem} 

\begin{theorem}\cite{K:86}
	Let $X_{j}$, $j=1,2,\ldots,n$, be any (not necessarily i.i.d.) random variables. For any $k=1,2,\ldots,n$ and $m=0,1,\ldots,n$, define random variables $Z_{m,k}^{n}$ such that:
	
	\item[(1)] $Z_{0,k}^{n}{\equiv}+\infty$, for $1{\leq}k{\leq}n$;
	\item[(2)] $Z_{m,k}^{n}{\equiv}-\infty$, for $m>n-k+1$;
	\item[(3)] $Z_{1,n}^{n}=X_{n}$;
	\item[(4)] $Z_{m,k}^{n}=[X_{k}{\vee}E[Z_{m,k+1}^{n}|F_{k}]]{\wedge}E[Z_{m-1,k+1}^{n}|F_{k}]$, for $1{\leq}m{\leq}n-k+1$, $k{\leq}n-1$.
	
	where $F_{k}$, $k=1,2,\ldots,n-1$, is a sigma-field over all possible realizations of vector $\{X_{i}\}_{i=1}^{k}$, $k=1,2,\ldots,n-1$, $\vee$ denotes the maximum, and $\wedge$ denotes the minimum.
\end{theorem}

If the performance rates of the workers are revealed at the beginning of the process then the problem reduces to the original case where the performance rates were known. A similar problem with i.i.d. random performance rates that follow a distribution $P$ at each stage of the problem, was studied for the SSAP in \cite{KHJ:16}. The main result is presented below

\begin{theorem}\cite{KHJ:16}\label{dssap_thm}
	The Greedy algorithm, which assigns the arriving task to the worker with the maximum performance rate value at each stage, achieves the maximum total expected reward in the DSSAP with i.i.d. random performance rates. Moreover, the maximum total expected reward is given by 
	
	\begin{equation}
	E[X]{\times}\big( \sum_{i=1}^{n} \int_{-\infty}^{+\infty} yd(F^{n-i+1}_{Q}(y)) \big)
	\end{equation}
	where $E[X]$ denotes the expected value of the random variables $X$ and $F^{n-i+1}_{Q}$ denotes the cdf of the workers' random performance rates at stage $i$ of the process where we have $n-i+1$ i.i.d. random performance rates.
\end{theorem}

\subsection{Threshold Reward Function} 
We introduce a new type of reward function (\ref{thres_cost_function}) defined using indicator functions which capture threshold constraints. Again, we consider $n$ jobs $\{X_{i},i=1,2,\ldots,n\}$ arriving sequentially, each following a distribution $G_{X_{i}}$, and $n$ workers with performance rates $p_{j}$, $ 0<p_{j}{\leq}1$. By the time job $i$ arrives, its value $x_{i}$ is randomly generated following $G_{X_{i}}$ and it is assigned online to a worker. The reward function is

\begin{equation}\label{thres_cost_function}
r(X)=\max\limits_{\pi}  [\sum_{i=1}^{n}\boldsymbol{1}_{\big\{\sum_{j=1}^{n}A_{ij}f(x_{i},p_{j}){\geq}\alpha \big\} }]
\end{equation}
the threshold $\alpha{\in}\mathbb{R}$, the threshold function $f$ is a two variable function $f:X{\times}P{\to}{\mathbb{R}}$, where $X=\{x_{1},x_{2}\ldots,x_{n}\}$ is the set of job values, $P=\{p_{1},p_{2}\ldots,p_{n}\}$ is the set of workers performance rates, and ${\pi}=\{A_{ij}{\in}\{0,1\}; \sum_{j=1}^{n}A_{ij}{\leq}1; i=1,2,\ldots,n; \sum_{i=1}^{n}A_{ij}{\leq}1; j=1,2,\ldots,n; \}$ represents the assignment policy. The objective is to find an optimal policy ${\pi}^{*}$ that determines the assignment of jobs to workers, $A_{ij}{\in}\{0,1\}$, $i=1,2,\ldots,n$, $j=1,2,\ldots,n$, such that $\sum_{j=1}^{n}A_{ij}{\leq}1$, $i=1,2,\ldots,n$, $\sum_{i=1}^{n}A_{ij}{\leq}1$, $j=1,2,\ldots,n$, and the reward function \eqref{thres_cost_function} is maximized.

For SSTAP we perform a stronger type of optimization. We do not maximize the expected value of the reward function but for each randomly generated sequence of job values we maximize the reward function. We focus our attention on a general class of threshold functions $f$, the order-preserving functions on the argument of worker performance rate.

\begin{definition} \textbf{(Order-preserving function)}
	Consider the two variable function $f:X{\times}P{\to}{\mathbb{R}}$, where $X=\{x_{1},\ldots,x_{n}\}$, $P=\{p_{1},\ldots,p_{n}\}$ are discrete sets of positive real numbers. The function $f$ is order-preserving on the arguments of $P$ if the order of the values $f(x_{i},p_{1})$, $f(x_{i},p_{2})$, $\ldots$, $f(x_{i},p_{n})$ is independent of the $x_{i}{\in}X$. 
\end{definition}

\section{Optimal Assignment Policy}

We present the optimal assignment policy ${\pi}^{*}$ given a threshold $\alpha$, job values $\{x_{1},x_{2},\ldots,x_{n}\}$ arriving sequentially and worker performance rates $p_{j}$, $j=1,2,\ldots,n$. The number of jobs, $n$ is fixed. The optimal policy is provided for order-preserving functions on the argument of worker value. Each time a job $x_{i}$ arrives, it is assigned to the non assigned worker $j$ with $p_{j}=arg\{\min\limits_{p}f(x_{i},p){\geq}\alpha\}$. If no such $p_{j}$ exists then $x_{i}$ is rejected. First, we give a lemma used in the proof of the theorem.

\begin{lemma}
	Let $p_{j_{k}}=arg\{\min\limits_{p}f(x_{k},p){\geq}{\alpha}\}$ then
	\begin{equation}\label{rel1}
	{\alpha}{\leq}f(x_{k},p_{j_{k}}){\leq}f(x_{k},p_{u})
	\end{equation} 
	\begin{equation}\label{rel2}
	f(x,p_{j_{k}}){\leq}f((x,p_{u}) \quad \forall x{\in}X 
	\end{equation}
	
	where $p_{u}{\neq}p_{j_k}$ is any performance rate such that $\alpha{\leq}f(x_{k},p_{u})$. 
\end{lemma}

\textbf{Proof:}
By definition, since $p_{j_{k}}=arg\{\min\limits_{p}f(x_{k},p){\geq}{\alpha}\}$ then for any $p_{u}$ such that $\alpha{\leq}f(x_{k},p_{u})$, it holds that ${\alpha}{\leq}f(x_{k},p_{j_{k}}){\leq}f(x_{k},p_{u})$. Since $f$ is order-preserving in the argument of the performance rates, we have that $f(x,p_{j_{k}}){\leq}f((x,p_{u})$ for any $x{\in}X$.
$\blacksquare$

Although the order-preserving function guarantees $f(x,p_{j_{k}}){\leq}f((x,p_{u})$ for any $x{\in}X$, it does not imply that ${\alpha}{\leq}f(x,p_{j_{k}})$ for any $x{\in}X$. We provide the main theorem of the paper.

\begin{theorem}\label{optimal_policy}
	Given a set of job values $X=\{x_{1},x_{2},\ldots,x_{n}\}$ arriving sequentially where $x_{i}$ is randomly generated following $G_{X_{i}}$, a set of performance rates $P=\{p_{1},p_{2},\ldots,p_{n}\}$, an order-preserving function $f$ and a threshold value $\alpha$, the optimal assignment policy ${\pi}^{*}$ that maximizes the cost function 
	\begin{equation}
	r(X)=\max\limits_{\pi}  [\sum_{i=1}^{n}\boldsymbol{1}_{\big\{\sum_{j=1}^{n}A_{ij}f(x_{i},p_{j}){\geq}\alpha \big\} }]
	\end{equation}
	is to assign to each arriving $x_{i}$ the not already assigned worker who corresponds to $p_{j}=arg\{\min\limits_{p}f(x_{i},p){\geq}\alpha\}$, where $p$ belongs to the set of the performance rates of the non-assigned workers. If no such $p_{j}$ exists then $x_{i}$ is rejected.
\end{theorem}

\textbf{Proof:}
Suppose that $m({\leq}n)$ jobs have arrived. Define the following sets: $X_{m}$ contains the first $m$ jobs that appeared, $S_{m}$ contains the performance rates of the already assigned workers, $N_{m}$ contains the performance rates of the non-assigned workers and $H_{m}{\subseteq}X_{m}$ contains the jobs assigned to some worker (not rejected). 

The proof proceeds by induction on $n$. For $n=1$, we assign the single job to the worker with $p_{j}=arg\min\limits_{p{\in}N_1}\{f(x_{i},p){\geq}\alpha\}$ and we get the maximum possible reward which is $1$. If such $p_{j}$ does not exist then the job is rejected and we get zero reward. Thus using  the suggested policy, we maximize the reward. We now assume that the claim holds for $n=k$, we prove that it holds for $n=k+1$.

Let $|H_{k}|=k'{\leq}k$, where $H_{k}=\{y_{1},y_{2},\ldots,y_{k'}\}$ are ordered and assigned to $k'$ workers with performance rates $\{p_{j_1},p_{j_2},\ldots,p_{j_{k'}}\}$, respectively. Let $x_{k+1}$ be the $(k+1)^{th}$ job value. If there is a performance rate $p_{i_{k+1}}{\in}N_{k}$ such that $p_{i_{k+1}}=arg\{\min\limits_{p{\in}N_{k}}{f(x_{k+1},p)}{\geq}\alpha\}$, then we assign the $(k+1)^{th}$ job to the $i_{k+1}$ worker.

Let us assume that such a success rate does not exist and that there is an alternate assignment of the $k$ preceding jobs such that the jobs $x_{k+1}{\cup}H_{k}$ can all pass the threshold. We recognize two assignments, the initial and the alternate assignment. In the initial assignment we follow the suggested policy and the jobs $\{y_{1},y_{2},\ldots,y_{k'}\}$ are assigned to the performance rates $\{p_{j_1},p_{j_2},\ldots,p_{j_{k'}}\}$, respectively. In the alternate assignment, the jobs $\{y_{1},y_{2},\ldots,y_{k'}\}$ are assigned to the performance rates $\{p_{g(j_{1})},p_{g(j_{2})},\ldots, p_{g(j_{k'})}\}$ respectively and $x_{k+1}$ to $p_{q}$; let $Q{\equiv}\{p_{g(j_{1})},p_{g(j_{2})},\ldots,p_{g(j_{k'})},p_{q}\}$, where $g$ is a mapping from $\{j_{1},j_{2},\ldots,j_{k'}\}$ to $\{1,2,\ldots,n\}$. 

If a job is successfully assigned to a worker in the initial assignment then it will also be assigned to a worker in the alternate assignment; not necessarily the same. This is because there exist suitable performance rates such that the order-preserving function evaluated at this job is greater than the threshold. It is meaningless to ignore a job that "passes" the threshold because we reduce the maximum reward by one unit that may not be replaced by one of the upcoming jobs. Even if it is replaced by one of the upcoming jobs we could have kept the job we initially ignored without changing the optimal reward. (PUT THIS AS A CONDITION IT THE PROBLEM STATEMENT).   

Since, in the alternate assignment we are able to assign job $x_{k+1}$ to a worker the reward $r_{alt}$, provided by the alternate assignment, will be increased by one compare to the reward $r_{init}$, provided by the initial assignment i.e., $r_{alt}=r_{init}+1$. We describe a process for all jobs $\{y_{1},y_{2},\ldots,y_{k'}\}$ according to which the performance rates $\{p_{g(j_{1})},p_{g(j_{2})},\ldots,p_{g(j_{k'})}\}$ of the alternate assignment can be swapped with the performance rates $\{p_{j_1},p_{j_2},\ldots,p_{j_{k'}}\}$ of the initial assignment without reducing the $r_{alt}$. We provide a detailed exposition of the process for job $y_{1}$ and it is the same for the jobs $y_{2},\ldots,y_{k'}$.

If $p_{j_{1}}\not\in Q$ we can assign $p_{j_{1}}$ to $y_{1}$ since it passes the threshold in the initial assignment, i.e $f(y_{1},p_{j_{1}}){\geq}{\alpha}$; and $p_{g(j_1)}$ is placed in $N_{k+1}$. If $p_{j_1}{\in}Q$ then there exists $l{\in}\{1,2,\ldots,k'\}$ such that $p_{j_{1}}=p_{g(j_{l})}$, where $p_{g(j_{l})}$ is assigned to job $y_{l}$ in the alternate assignment. From the alternate assignment we have $f(y_{1},p_{g(j_1)}){\geq}\alpha$ and $f(y_{l},p_{g(j_l)}){\geq}\alpha$. We claim that we can swap the success rates $p_{g(j_{1})}$, $p_{g(j_{l})}=p_{j_{1}}$ and still the thresholds are satisfied i.e., $f(y_{1},p_{g(j_l)}){\geq}\alpha$ and $f(y_{l},p_{g(j_1)}){\geq}\alpha$.

In the initial assignment, we apply the suggested policy which implies: 
\begin{equation}\label{rel3}
p_{j_1}=arg\{\min\limits_{p}f(y_{1},p){\geq}\alpha\}
\end{equation}
For job $y_{1}$, the \eqref{rel3} and \eqref{rel1} imply
\begin{equation}\label{swap1}
f(y_{1},p_{j_1}=p_{g(j_l)}){\geq}\alpha
\end{equation}
For job $y_l$, from the alternate assignment we have
\begin{equation}\label{rel4}
f(y_{l},p_{g(j_l)}){\geq}\alpha
\end{equation}
The \eqref{rel3} and \eqref{rel1} imply $f(y_1,p_{j_1}){\leq}f(y_1,p_{g(j_1)})$ and \eqref{rel2} tells that $f(x,p_{j_{1}}=p_{g(j_l)}){\leq}f((x,p_{g(j_1)}) \quad \forall x{\in}X$. For $x=y_l$, $f(y_l,p_{j_1}=p_{g(j_l)}){\leq}f(y_l,p_{g(j_1)})$ and from \eqref{rel4}, we have 
\begin{equation}\label{swap2}
f(y_l,p_{g(j_1)}){\geq}\alpha
\end{equation}
This is the end of the process. Equations \eqref{swap1}, \eqref{swap2} prove the claim. We continue this process for the success rates of all jobs: $y_{2}$,$\ldots$,$y_{k'}$. When we finish with all jobs, we get the following assignment $\{p_{j_{1}}, p_{j_{2}},\ldots, p_{j_{k'}}\}$, respectively and $x_{k+1}$ to $p_{q'}$, where $p_{q'}$ is not necessarily equal to $p_{q}$, due to the swaps that take place in the process, in any case $f(x_{k+1},p_{q'}){\geq}\alpha$. Therefore $p_{q'}\in S_{k+1}$, which is a contradiction, since we assumed that $p_{q'}=arg\{\min\limits_{p{\in}N_{k}}{f(x_{k+1},p)}{\geq}\alpha\}$ does not exist. Hence, it is optimal the $x_{k+1}$ to be assigned to a performance rate according to the suggested policy. The claim holds for $n=k+1$. This concludes the induction step. The suggested policy is optimal.
$\blacksquare$

We provide two examples. The first one is an application of the optimal policy. The second example shows the necessity of the order-preserving function in order the suggested policy to be optimal. 

\begin{example}
	Consider the SSTAP for four workers and four jobs with the order-preserving function $f(x,p)=xp$, threshold $\alpha=0.15$ and performance rates $p_{1}=0.4$, $p_{2}=0.5$, $p_{3}=0.6$, $p_{4}=0.7$ that imply the following ordering $f(.,p_{1}){\leq}f(.,p_{2}){\leq}f(.,p_{3}){\leq}f(.,p_{4})$. For the job values $x_{1}=0.0975$, $x_{2}=0.275$, $x_{3}=0.9575$, $x_{4}=0.4854$, arriving sequentially, we get the following assignment: $x_1$ rejected, $x_{2}$ assigned to $p_{3}$, $x_{3}$ assigned to $p_{1}$ and $x_{4}$ assigned to $p_{2}$.
\end{example}

\begin{example}
	We provide an example that highlights the necessity of an order-preserving function in SSTAP. For the performance rates $p_{1}$, $p_{2}$, $p_{3}$ and the job values $x_{1}$, $x_{2}$, $x_{3}$ arriving sequentially in this order, we have the following function in the threshold constraints:\\
	
	\begin{center}
		\begin{tabular}{ |c|c|c| } 
			\hline
			$f(x_{1},p_{1})=0.5$ & $f(x_{2},p_{1})=0.08$ & $f(x_{3},p_{1})=0.5$ \\ 
			\hline
			$f(x_{1},p_{2})=0.4$ & $f(x_{2},p_{2})=0.1$ & $f(x_{3},p_{2})=0.4$ \\ 
			\hline
			$f(x_{1},p_{3})=0.7$ & $f(x_{2},p_{3})=0.03$ & $f(x_{3},p_{3})=0.1$ \\ 
			\hline
		\end{tabular}
	\end{center}
	\vspace{0.5cm}
	The $f$ is not order-preserving, since for $x_{1}$ $f(x_{1},p_{2}){\leq}f(x_{1},p_{1}){\leq}f(x_{1},p_{3})$ while for $x_{2}$ $f(x_{2},p_{3}){\leq}f(x_{2},p_{1}){\leq}f(x_{2},p_{2})$. For threshold $\alpha=0.1$ we observe that according to the suggested policy $x_{1}$ is assigned to $p_{2}$. The $x_{2}$ is aborted since only $p_{2}$ could satisfy the threshold, but it is already assigned. The $x_{3}$ is assigned to $p_{3}$. The reward is $2$. However, we could have assigned $x_{1}$ to $p_{1}$, $x_{2}$ to $p_{2}$ and $x_{3}$ to $p_{3}$ and get a reward of $3$.
\end{example}

\subsection{Illustrative Example}

We give an example of the optimal assignment policy for a reward function with threshold constraints inspired from aviation security applications. We assume that the job values $x_{i}{\in}[0,1]$ stand for the risk value of the passenger and the performance rates $p_{i}{\in}[0,1]$ quantifies the capabilities of the workers, for $i=1,2,\ldots,n$. A job $x{\approx}0$ stands for a low risk passenger and he must be assigned to an officer of lower capabilities $p{\approx}0$. Similarly, a job $x{\approx}1$ represents a high risk passenger and he must be assigned to high capabilities officer, $p{\approx}1$. To this end, we introduce the following cost function,

\begin{equation}
r_{av}=max[\sum_{i=1}^{n}\boldsymbol{1}_{\big\{\sum_{j=1}^{n}\frac{1}{x_{i}}A_{ij}p_{j}{\geq}a \big\} }]
\end{equation}

We observe that the threshold function $f(x_{i},p_{j})=\frac{p_{j}}{x_{i}}$ is order-preserving in the argument of $p's$, thus we can apply the optimal policy algorithm. For a threshold value $\alpha$ in $\frac{p_{j}}{x_{i}}{\geq}{\alpha}$, if $x_{i}{\approx}0$ the job passes the threshold with very high probability and it will be served by a low performance rate officer, which is what we expect for a low risk passenger. On the other hand, if $x_{i}{\approx}1$ we need a higher performance rate officer in order to pass the threshold, which also describes the problem appropriately.

We provide a figure which depicts the number of passengers which pass the threshold out of a total of 200 passengers. The threshold varies from 0.1 to 5. The values of the jobs follow the uniform distribution $U(0,1)$. The values of the 200 workers are given by the expression $p_{i}=\frac{i}{200}$, for $i=1,\ldots,200$.
\begin{figure}
	\centering
	\includegraphics[scale=0.6]{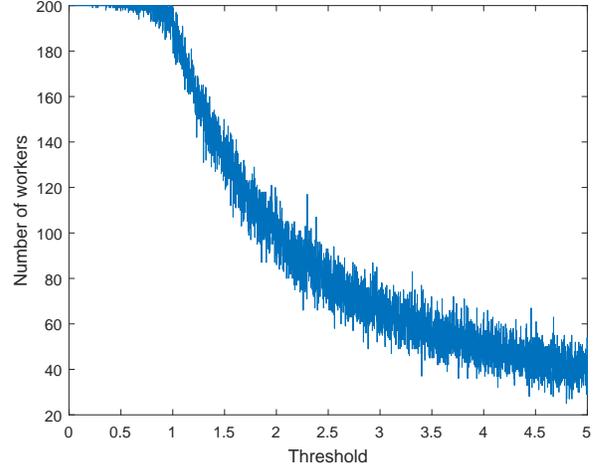}
	\caption{Number of passengers which pass the threshold out of a total of 200 passengers, as the threshold varies from 0.1 to 5.}
\end{figure}
For threshold $\alpha=0.5$, the jobs which pass the threshold are 200  out of 200. For threshold $\alpha=1$, the jobs which pass the threshold are 193 out of 200. For threshold $\alpha=3$, the jobs which pass the threshold are 58 out of 200.

\subsection{Infinite Number of Jobs and Workers Cycling Back}

In passenger screening systems for aviation security, the total number of passengers that will be examined is a random variable. In Internet transaction, the sequentially arriving tasks are infinite and the available servers must handle them. The optimal policy works for random number of jobs, even infinite as long as there exist available workers. The lack of available workers can be resolved by accepting workers that cycle back. 

In cycling back, workers can be reused after the completion of a task. Each worker cycles back with cycle rate ${\lambda}_{i}$ for $i=1,2,\ldots,M$, where $M$ is the number of workers. By definition, if worker $i$ never cycles back, then ${\lambda}_{i}=\infty$. The set of available workers is updated when a worker returns. We apply the policy given by Theorem \eqref{optimal_policy} at the arrival of a new job and we use as input the updated list of available workers.   

\section{Performance analysis}\label{sec_perf_an}
In this section, we provide a performance analysis of a system that uses SSTAP optimal assignment policy. We research the maximum and minimum \textit{job load} that a SSTAP system can service while achieving its maximum reward. The maximum and minimum job load are marginal values which if exceeded or missed respectively, the maximum reward is reduced. 
\begin{definition}\label{job_load}
	We define the job load, for the set of job values $X=\{x_{1},x_{2},\ldots,x_{n}\}$, as the Euclidean norm
	\begin{equation}
	l(X)=\sqrt{x_{1}^{2}+x_{2}^{2}+\ldots+x_{n}^{2}}
	\end{equation}
\end{definition}

Given the performance rates $\{p_{i},i=1,2,\ldots,n\}$ and an order-preserving function $f$, we compute the set $M=\{u_{i}=arg\max\limits_{x}\{\alpha{\leq}f(x,p_{i})\},x{\in}\Omega,i=1,2,\ldots,n\}$. The $l(M)$ is the maximum job load for reward equal to $n$. If we increase the job load, the reward will be reduced. By computing the set $N=\{v_{i}=arg\min\limits_{x}\{\alpha{\leq}f(x,p_{i})\},x{\in}\Omega,i=1,2,\ldots,n\}$, the $l(N)$ is the minimum job load for reward equal to $n$. If we further decrease the job load, the reward will be reduced.

\begin{theorem}
	We consider the set of $n$ job values $X=\{x_{1},x_{2},\ldots,x_{n}\}$, a threshold value $\alpha$, an order-preserving function $f$ and the sets $M=\{u_{i}=arg\max\limits_{x}\{\alpha{\leq}f(x,p_{i})\},x{\in}\Omega,i=1,2,\ldots,n\}$, $N=\{v_{i}=arg\min\limits_{x}\{\alpha{\leq}f(x,p_{i})\},x{\in}\Omega,i=1,2,\ldots,n\}$. Then
	
	\begin{equation}
	r(X)=n {\Rightarrow} l(N){\leq}l(X){\leq}l(M)
	\end{equation}
	where $l$ is the job load defined in Definition \ref{job_load} and $r(X)$ is the value of the reward function \eqref{thres_cost_function} given that for the sequential assignment of the jobs $X=\{x_{1},x_{2},\ldots,x_{n}\}$ we follow the optimal assignment policy given in Theorem \ref{optimal_policy}. 
\end{theorem}

\textbf{Proof:}
We proceed by contradiction. Let us assume that $r(X)=n$, this implies that all jobs with values $X=\{x_{1},x_{2},\ldots,x_{n}\}$ are assigned to workers with performance rates $\{p_{i_1},p_{i_2},\ldots,p_{i_n}\}$ respectively.

a) If we assume that $l(X)>l(M)$ then there exists $x_{j}{\in}X$ such that $x_{j}>u_{i_j}=arg\max\limits_{x}\{\alpha{\leq}f(x,p_{i_j})\}$, $x{\in}\Omega$ $\Rightarrow$ $f(x_j,p_{i_j})<\alpha$ $\Rightarrow$ $r(X)<n$, which is a contradiction.

b) If we assume that $l(X)<l(N)$ then there exists $x_{j}{\in}X$ such that $x_{j}<v_{i_j}=arg\min\limits_{x}\{\alpha{\leq}f(x_j,p_{i_j})\}$, $x{\in}\Omega$ $\Rightarrow$ $f(x_j,p_{i_j})<\alpha$ $\Rightarrow$ $r(X)<n$, which is a contradiction.
$\blacksquare$

In the SSAP the random job values $\{X_{i},i=1,2,\ldots,n\}$ that follow $G_{X_{i}}$ determine the subintervals in the domain $\Omega$ of $\{X_{i},i=1,2,\ldots,n\}$ for the optimal sequential assignment. In SSTAP the optimal sequential assignment algorithm does not depend on the distributions of the job values. Taking advantage of that, we can look for a probability distribution $G_{X}$ for the i.i.d. random job values $\{X_{i},i=1,2,\ldots,n\}$ that maximizes the reward function \eqref{thres_cost_function} for the maximum and minimum job load respectively.

We take $u_{i}'=u_{i}-{\epsilon}$, for $\epsilon$ small positive number, the mass probability function $P(x=u_{i}')=\frac{1}{n}$, $i=1,2,\ldots,n$ maximizes the reward function given that we use the optimal policy. In case we have $k$ equal values $u_{i_{1}}'=\ldots=u_{i_{k}}'$, $P(x=u_{i_{1}}')=\frac{k}{n}$. However, we are interested in continuous probability distributions. To this end, we approximate the mass probability function using Gaussian distributions $N_{i}$ with mean values $u_{i}$ and standard deviations very small ${\sigma}{\approx}0$ and peak value at $P(x=u_{i}')$

\begin{equation}\label{gaussian}
f_{\sigma}(x)=\frac{1}{w}\sum_{i}N_{i}(x_{i},\sigma)(x)
\end{equation}
where $w=\int_{x{\in}\Omega}\sum_{i}N_{i}(x_{i},\sigma)(x)dx$. Based on the definition \eqref{conv_dist}, we observe that as ${\sigma}{\to}0$ the distribution $f(x)$ converges in distribution to the mass probability function which maximizes the reward function. Similar results to \eqref{gaussian} could be reached for the set $N$ by taking $v_{i}'=v_{i}+{\epsilon}$,

\begin{figure}\label{figu1}
	\centering
	\includegraphics[scale=0.6]{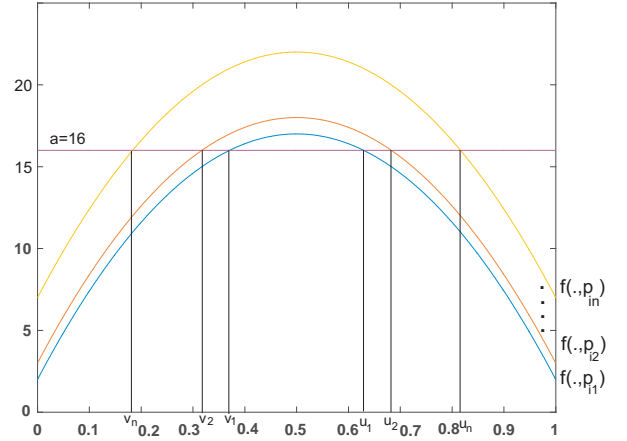}
	\caption{Intervals we use for the choice of the $x_{i}'s$.}
\end{figure}

\section{Multiple Levels SSTAP}
The screening process of the airplane passengers is performed in multiple levels. This fact inspired the study of the multiple levels SSTAP. Multiple levels SSTAP is a generalization of SSTAP, where workers are organized into multiple levels. For SSTAP with $k$-levels, a subset of workers is operating in a single level. In the $i^{th}$ level, $i=1,2,\ldots,k$, there are $s_{i}$ workers with performance rates $\{p_{j_{1}},p_{j_{2}},\ldots,p_{j_{s_{i}}}\}$, threshold ${\alpha}_{i}$ and the order-preserving function $f_{i}$. The reward function in the $i^{th}$ level is given by

\begin{equation}
r_{i}=\max\limits_{{\pi}_{i}}  [\sum_{j=1}^{s_{i}}\boldsymbol{1}_{\big\{\sum_{w=1}^{s_{i}}A_{jw}f_{i}(x_{j},p_{w}){\geq}{\alpha}_{i} \big\} }]
\end{equation}

For an arriving job $x$, the optimal policy can be applied at level $1$. If the job is rejected in level $1$, it proceeds to level $2$. The process continues until we reach the last level, where if the job is rejected, it is never assigned. The multiple levels SSTAP is characterized by a priority property from lower to higher levels. If a job can be assigned to a worker of the $i^{th}$ level then it must be assigned without proceeding to the next levels. The goal of SSTAP is to maximize the $\sum_{i=1}^{n}r_{i}$ under the priority property. We provide the theorem that describes the optimal assignment policy for the $k$-level SSTAP.

\begin{theorem}
	The optimal assignment policy for the $k$-level SSTAP, is the vector $\tilde{\pi}^{*}=({\pi}^{*}_{1},{\pi}^{*}_{2},\ldots,{\pi}^{*}_{k})$, where ${\pi}^{*}_{i}$ is the optimal policy given by the algorithm described in Theorem \ref{optimal_policy} for the $i^{th}$ level of the problem, $1{\leq}i{\leq}k$.
\end{theorem}

\textbf{Proof:}
With respect to the priority property, we apply the policy described in Theorem \ref{optimal_policy} at each level and we get the vector of the optimal assignment policy for the $k$-levels problem $\mathbf{{\pi}^{*}}=({\pi}^{*}_{1},{\pi}^{*}_{2},\ldots,{\pi}^{*}_{k})$. If a job fails to be assigned at its current level it moves to the next level. At the final level if a job fails to be assigned, it is aborted.
$\blacksquare$

The following lemma states that the partition of workers into multiple levels under the priority property may result in a smaller reward compared to the single level case.

\begin{lemma}
	Consider the $k$-levels SSTAP, $k{\geq}2$, with the same order-preserving function in every level, $f_{i}=f$, and the same threshold ${\alpha}_{i}=\alpha$, $1{\leq}i{\leq}k$. We consider the induced $1$-level SSTAP with reward function
	\begin{equation}
	r=\sum_{i=1}^{k}\sum_{j=1}^{s_{i}}\boldsymbol{1}_{\big\{\sum_{w=1}^{s_{i}}A_{jw}f(x_{j},p_{w}){\geq}{\alpha} \big\} }
	\end{equation}
	For the total reward $\sum_{i=1}^{k}r_{i}$ of the $k$-levels SSTAP under the optimal policy $\tilde{\pi}^{*}$ and the reward $r$ of the induced $1$-level SSTAP under the optimal policy ${\pi}^{*}$, it holds
	\begin{equation}
	\sum_{i=1}^{k}r_{i}{\leq}r
	\end{equation}
\end{lemma}

\textbf{Proof:}
\begin{equation}
\sum_{i=1}^{k}r_{i}\Big|_{\tilde{\pi}^{*}}=\sum_{i=1}^{k}\sum_{j=1}^{s_{i}}\boldsymbol{1}_{\big\{\sum_{w=1}^{s_{i}}A_{jw}f(x_{j},p_{w}){\geq}{\alpha} \big\}} \Big|_{\tilde{\pi}^{*}}{\leq}\sum_{i=1}^{k}\sum_{j=1}^{s_{i}}\boldsymbol{1}_{\big\{\sum_{w=1}^{s_{i}}A_{jw}f(x_{j},p_{w}){\geq}{\alpha} \big\}} \Big|_{{\pi}^{*}}
\end{equation}
where the last inequality comes from the optimality of ${\pi}^{*}$.
$\blacksquare$

The multiple levels SSTAP permit us to organize the workers into groups according to their performance rates. For example, for a three level multiple levels SSTAP, following the ordering indicated by the order-preserving function $f(.,p_{i_1}){\leq}f(.,p_{i_{2}}){\leq}{\ldots}{\leq}f(.,p_{i_{n}})$, we can place the first $70\%$ of the workers in the first level, the next $20\%$ to the second level and the remaining $10\%$ to the third level.

\section{SSTAP with Random Performance Rates}\label{sec_sstap_rpr}
We extend the result \eqref{dssap_thm} for SSTAP. This new problem is denoted as doubly stochastic sequential threshold assignment problem (DSSTAP). The reward function, we maximize is

\begin{equation}\label{rewa_perf_rates}
\begin{split}
&\max\limits_{\pi}E[ \sum_{i=1}^{n}\boldsymbol{1}_{\big\{\sum_{j=1}^{n}A_{ij}f(x_{i},p_{j}){\geq}\alpha \big\} }]\\
&=\max\limits_{\pi}\sum_{i=1}^{n}Pr(\sum_{i=1}^{n}A_{ij}f(x_{i},p_{j}){\geq}\alpha)
\end{split}
\end{equation}

\textbf{Case I:} We consider i.i.d. job values that follow the distribution $X{\sim}G_{X}$ and worker performance rates not i.i.d., that follow the distribution $P_{i}{\sim}G_{P_{i}}$. The reward function is independent of the policy $\pi$ we apply

\begin{equation}
\begin{split}
&\max\limits_{\pi}E[ \sum_{i=1}^{n}\boldsymbol{1}_{\big\{\sum_{j=1}^{n}A_{ij}f(x_{i},p_{j}){\geq}\alpha \big\} }]\\
&=\sum_{i=1}^{n}Pr(f(X,P_{i}){\geq}\alpha)
\end{split}
\end{equation}

\textbf{Case II:} We consider job values not i.i.d. that follow the distribution $X_{i}{\sim}G_{X_{i}}$ and worker performance rates not i.i.d., that follow the distribution $P_{i}{\sim}G_{P_{i}}$. The reward function depends on the policy we apply and the assignment problem reduces to the maximum weighted matching of the bipartite graph between the disjoint sets $X=\{x_{1},x_{2},\ldots,x_{n}\}$ and $P=\{p_{1},p_{2},\ldots,p_{n}\}$ with edge weights $w_{ij}=Pr(f(X_{i},P_{j}){\geq}\alpha)$. The maximum weighted matching of the bipartite graph is resolved using the Hungarian algorithm, in $O(n^4)$ time complexity.

\section{Conclusion}
We introduce SSTAP, a variation of the SSAP problem defined using indicator functions that capture threshold constraints. An optimal assignment policy is proven and it is independent from the number of jobs. Via an illustrative example, we show that SSTAP models accurately aviation security problems. We provide a performance analysis of systems that use SSTAP optimal policy for their sequential assignment projects. Finally, we analyze the multiple levels SSTAP, and we study the SSTAP with random performance rates.\\


\textbf{Acknowledgement}:I would like to thank professor Sheldon Howard Jacobson (Department of Computer Science, University of Illinois at Urbana Champaign) for the insightful disucssions and for providing me research questions presented in this paper.

\end{document}